\documentclass[a4paper,11pt]{article}
\usepackage{amsfonts}
\usepackage{mathbbold}%
\usepackage{amsfonts}
\usepackage{bbding}
\usepackage{amssymb}
\usepackage{mathrsfs}
\usepackage{float}
\usepackage{amsmath}
\usepackage{amsfonts,amssymb}
\usepackage{graphicx,color}
\usepackage{psfig}
\usepackage{mathrsfs}
\usepackage{graphicx}%
\catcode`\@=11
\@addtoreset{equation}{section}

\catcode`\@=12

\newtheorem{Theorem}{Theorem}[section]

\newtheorem{Lemma}{Lemma}[section]

\newtheorem{Remark}{Remark}[section]

\newtheorem{Proposition}{Proposition}[section]

\def\2{{I \hskip -1.0mm I}}
\def\3{{I \hskip -1.0mm I\hskip -1.0mm I}}
\def\4{{I \hskip -0.9mm V}}
\def\6{{V \hskip -1.35mm I}}

\begin{document}
\title{Life-span of classical solutions to hyperbolic geometric flow in two space variables
with slow decay initial data}
\author{De-Xing Kong$^{a}$, $\quad$ Kefeng Liu$^{b}$ $\quad$ and $\quad$ Yu-Zhu
Wang$^{c}$%\footnote{Corresponding author.}
\\ \\
$^a\!${\it {\small Department of Mathematics, Zhejiang University}}\\
{\it  {\small Hangzhou 310027, China}}\\
$^b\!${\it {\small Department of Mathematics, University of
California at Los
Angeles}}\\
{\it  {\small CA 90095, USA}}
\\$^{c}\!${\it {\small
Department of Mathematics, Shanghai Jiao Tong University}}\\{\it
{\small Shanghai 200248, China}}}
\date{}
\maketitle

\begin{abstract}
In this paper we investigate the life-span of classical solutions to
the hyperbolic geometric flow in two space variables with slow decay
initial data. By establishing some new estimates on the solutions of
linear wave equations in two space variables, we give a lower bound
of the life-span of classical solutions to the hyperbolic geometric
flow with asymptotic flat initial Riemann surfaces.

\vskip 6mm


\vskip 3mm


\end{abstract}

\newpage\baselineskip=7mm

\section{Introduction}
Let $\mathscr{M}$  be an $n$-dimensional complete Riemannian
manifold with Riemannian metric $g_{ij}$. The following evolutionary
equation for the metric $g_{ij}$
\begin{equation}
\label{1.1}\dfrac{\partial^{2}g_{ij}}{\partial t^{2}}=-2R_{ij}
\end{equation}
has been recently introduced by Kong and Liu \cite{kl} and named as
{\it hyperbolic geometric flow}, where $R_{ij}$ stands for the Ricci
curvature tensor of $g_{ij}$. For the study on the hyperbolic
geometric flow, we refer to the recent papers \cite{dkl1},
\cite{dkl2}, \cite{k2}, \cite{kl} and \cite{klx}.

We are interested in the evolution of a Riemannian metric $g_{ij}$
on a Riemann surface $\mathscr{S}$ under the flow (\ref{1.1}). On a
surface, the hyperbolic geometric flow equation (\ref{1.1})
simplifies, because all of the information about curvature is
contained in the scalar curvature function $R.$ In our notation,
$R=2K$ where $K$ is the Gauss curvature. The Ricci curvature is
given by
\begin{equation}
\label{1.2}R_{ij}=\dfrac12Rg_{ij},
\end{equation}
and the hyperbolic geometric flow equation (\ref{1.1}) simplifies
the following equation for the special metric
\begin{equation}
\label{1.3}\dfrac{\partial^{2}g_{ij}}{\partial t^{2}}=-Rg_{ij}.
\end{equation}
The metric for a surface can always be written (at least locally) in
the following form
\begin{equation}
\label{1.4}g_{ij}=v(t,x,y)\delta_{ij},
\end{equation}
where $v(t,x,y)>0.$ Therefore, we have
\begin{equation}
\label{1.5}R=-\dfrac{\triangle\ln v}{v}.
\end{equation}
Thus the equation (\ref{1.3}) becomes
\[
\dfrac{\partial^{2}v}{\partial t^{2}}=\dfrac{\triangle\ln v}{v}\cdot
v,
\]
namely,
\begin{equation}
\label{1.6}v_{tt}-\triangle\ln v=0.
\end{equation}
Denote
\begin{equation}
\label{1.7}u=\ln v,
\end{equation}
then the wave equation (\ref{1.6}) reduces to
\begin{equation}
\label{1.8} u_{tt}-e^{-u}\Delta u=-u^2_t.
\end{equation}
(\ref{1.8}) is a quasilinear hyperbolic wave equation. The global
existence and the life-span of classical solutions to the Cauchy
problem for  hyperbolic equations  with the initial data with
compact support have been studied by many authors (e.g., \cite{k1},
\cite{lz}, \cite{li}, etc.). However, only a few results have been
known for the case of the initial data with non-compact support,
which plays an important role in both mathematics and physics.

Recently, Kong, Liu and Xu \cite{klx} studies the evolution of a
Riemannian metric $g_{ij}$ on a cylinder $\mathscr{C}$ under the
hyperbolic geometric flow (\ref{1.1}). They prove that, for any
given initial metric on ${\mathbb{R}}^{2}$ in a class of cylinder
metrics, one can always choose suitable initial velocity symmetric
tensor such that the solution exists for all time, and the scalar
curvature corresponding to the solution metric $g_{ij}$ keeps
uniformly bounded for all time; moreover, if the initial velocity
tensor is suitably ``large", then the solution metric $g_{ij}$
converges to the flat metric at an algebraic rate. If the initial
velocity tensor does not satisfy the condition, then the solution
blows up at a finite time, and the scalar curvature $R(t,x)$ goes to
positive infinity as $(t,x)$ tends to the blowup points, and a flow
with surgery has to be considered. This result shows that, by
comparing to Ricci flow, the hyperbolic geometric flow has the
following advantage: the surgery technique may be replaced by
choosing suitable initial velocity tensor. Some geometric properties
of hyperbolic geometric flow on general open and closed Riemann
surfaces are also discussed (see Kong et al \cite{klx}).

In this paper, we consider the Cauchy problem for (\ref{1.8}) with
the following initial data
\begin{equation}
\label{1.9} t=0:\;\;  u=\varepsilon u_0(x),\quad u_t=\varepsilon
u_1(x),
\end{equation}
where $\varepsilon > 0$ is a suitably small parameter, $u_0(x)$ and
$u_1(x)$ are two smooth functions of $x\in \mathbb{R}^2$ and satisfy
that there exist two positive constants $A\in \mathbb{R}^+$ and
$k>1, k\in \mathbb{R}^{+}$ such that
\begin{equation}
\label{1.10} |u_0(x)|\leq \frac{A}{(1+|x|)^k},\quad |u_1(x)|\leq
\frac{A}{(1+|x|)^{k+1}}.
\end{equation}
(\ref{1.10}) implies that the initial data satisfies the slow decay
property, that is, the initial Riemann surface are asymptotic flat.
We shall prove the following theorem.
\begin{Theorem}
Suppose that $u_0(x),\; u_1(x)\in C^\infty(\mathbb{R}^2)$ and
satisfy the decay condition (\ref{1.10}). Then there exist two
positive constants $\delta $ and $\varepsilon_0$ such that for any
fixed $ \varepsilon\in [0,\varepsilon_0]$, the Cauchy problem
(\ref{1.8})-(\ref{1.9}) has a unique $C^\infty$ solution on the
interval $[0, T_\varepsilon]$, where $T_\varepsilon$ is given by
\begin{equation}\label{1.11}
T_\varepsilon=\frac{\delta}{\varepsilon^{\frac{4}{3}}}.
\end{equation}
\end{Theorem}

As we know, the flow equation (\ref{1.1}) is a system of fully
nonlinear partial differential equations of second order, it is very
difficult to study the global existence or blow-up of the classical
solutions of (\ref{1.1}). An interesting and important question is
to investigate the evolution of asymptotic flat initial Riemann
surfaces under the flow (\ref{1.1}). In this case, although the
equation (\ref{1.1}) can simply reduce to (\ref{1.8}), (\ref{1.8})
is still a fully nonlinear wave equation, only a few results have
been known even for its Cauchy problem. Our main result, Theorem
1.1, gives a lower bound on the life-span of the classical solution
of the Cauchy problem (\ref{1.8})-(\ref{1.9}). This theorem shows
that the smooth evolution of asymptotic flat initial Riemann
surfaces under the flow (\ref{1.1}) exists at least on the interval
$[0, T_\varepsilon]$.

The paper is organized as follows. In Section 2 we establish some
new estimates on the solutions of linear wave equations in two space
variables, these estimates play an important role in the proof of
Theorem 1.1. Based on this, we prove Theorem 1.1 in Section 3, which
gives a lower bound of the life-span of classical solutions to the
hyperbolic geometric flow with asymptotic flat initial Riemann
surfaces.

\section{Some useful lemmas}
Following Klainerman \cite{kl2}, we introduce a set of partial
differential operators
\begin{equation}
Z=\{\partial_i\;\; (i=0, 1, \cdots, n);\quad L_{0};\quad
\Omega_{ij}\;\; (1\leq i<j\leq n); \quad \Omega_{0i}\;\; (i=1,
\cdots, n)\},
\end{equation}
where
\begin{equation}
\partial_0=\frac{\partial}{\partial t}, \quad \partial_i=\frac{\partial}{\partial
x_i}\quad (i=1, \cdots, n),
\end{equation}
\begin{equation}
L_0 =t\partial_0+\sum^n_{i=1}x_i \partial_i,
\end{equation}
\begin{equation}
\Omega_{ij}=x_i\partial_j-x_j\partial_i \quad ( 1\leq i<j\leq n)
\end{equation}
and
\begin{equation}
\Omega_{0i}=t\partial_i+x_i\partial_0 \quad (i=1, \cdots, n ).
\end{equation}
Let $Z^I$ denote a product of $|I|$ of the vector fields
(2.2)-(2.5), where $I=(I_1,\cdots, I_\sigma)$ is a multi-index,
$|I|=I_1+\cdots +I_\sigma$, $\sigma $ is the number of partial
differential operators in $Z: \; Z=(Z_1, \cdots, Z_\sigma)$ and
\begin{equation}
Z^I=Z^{I_1}_1\cdots Z^{I_\sigma}_\sigma.
\end{equation}

Throughout this paper, we use the following notations:
$L^p(\mathbb{R}^n)\; (1\leq p\leq \infty)$ stands for the usual
space of all $L^p(\mathbb{R}^n)$ functions on $ \mathbb{R}^n$ with
the norm $\|f\|_{L^p}$, $H^s$ denotes $s$-order Sobolev space on
$\mathbb{R}^n$ with the norm
$$\|f\|_{H^s}=\|(1+|\xi|)^{\frac{s}{2}}\hat{f}\|_{L^2},$$ where $s$ is a
given real number.

The following lemma has been proved in Li and Zhou \cite {lz}.

\begin{Lemma}
For any given multi-index $I=(I_1,\cdots, I_\sigma),$ we have
\begin{equation}
[\square, Z^I]=\sum_{|J|\leq|I|-1}A_{IJ}Z^J\square
\end{equation}
and
\begin{equation}
[\partial_{i},
Z^I]=\sum_{|J|\leq|I|-1}B_{IJ}Z^J\partial=\sum_{|J|\leq|I|-1}\tilde{B}_{IJ}\partial
Z^J\quad (i=0, 1, \cdots, n),
\end{equation}
where $ [\cdot ,\cdot]$ stands for the Poisson bracket, $J=(J_1,
\cdots, J_\sigma)$ a multi-index, $\square$ denotes the wave
operator, ${\displaystyle\partial=\left(\frac{\partial}{\partial t},
\frac{\partial}{\partial x_1}, \cdots, \frac{\partial}{\partial
x_n}\right)}$ and $ A_{IJ}, B_{IJ}$, $\tilde{B}_{IJ} $ stand for
constants.
\end{Lemma}

\begin{Lemma}
Assume that $n\geq 1.$ Let $u$ be a solution of the following Cauchy
problem
\begin{equation}\left\{\begin{array}{l}
\phi_{tt}-\triangle \phi=f,\\
t=0: \;\; u=\phi_0(x), \quad u_t=\phi_1(x).\end{array}\right.
\end{equation}
Then
\begin{equation} \|\partial \phi(t, \cdot)\|_{H^s}\leq
C(\|\partial_x \phi_0\|_{H^s}+\|\phi_1\|_{H^s}+\int^t_0\|f(\tau,
\cdot)\|_{H^s}),
\end{equation}
provided that all norms appearing in the right-hand side of (2.10)
are bounded.
\end{Lemma}

\noindent {\bf Proof.} Taking the Fourier transformation on the
variable $ x $ in (2.9) leads to
\begin{equation}
\left\{\begin{array}{l}
\hat{\phi}_{tt}+|\xi|^2\hat{\phi}=\hat{f}(t, \xi),\\
t=0: \hat{\phi}=\hat{\phi}_0(\xi),
\hat{\phi}_t=\hat{\phi}_1(\xi).\end{array}\right.
\end{equation}
Solving the initial value problem (2.11) gives
\begin{equation}
\hat{\phi}(t,
\xi)=\cos(t|\xi|)\hat{\phi}_0(\xi)+\frac{\sin(t|\xi|)}{|\xi|}\hat{\phi}_1(\xi)+\int^t_0\frac{\sin((t-\tau)|\xi|)}{|\xi|}\hat{f}(\tau,
\xi)d\tau.
\end{equation}
Thanks to (2.12), we obtain
\begin{equation}
\partial_t\hat{\phi}(t,
\xi)=-|\xi|\sin(t|\xi|)\hat{\phi}_0(\xi)+\cos(t|\xi|)\hat{\phi}_1(\xi)+\int^t_0\cos((t-\tau)|\xi|)\hat{f}(\tau,
\xi)d\tau
\end{equation}
and
\begin{equation}
|\xi|\hat{\phi}(t,
\xi)=|\xi|\cos(t|\xi|)\hat{\phi}_0(\xi)+\sin(t|\xi|)\hat{\phi}_1(\xi)+\int^t_0\sin((t-\tau)|\xi|)\hat{f}(\tau,
\xi)d\tau.
\end{equation}
It follows from  (2.13) and Minkowski inequality that
\begin{equation}\begin{array}{lll}
\|\partial_t\phi(t, \cdot)\|_{H^s} & \leq & {\displaystyle
\|(1+|\xi|)^{\frac{s}{2}}|\xi|\sin(t|\xi|)\hat{\phi}_0(\xi)\|_{L^2}+\|(1+|\xi|)^{\frac{s}{2}}\cos(t|\xi|)
\hat{\phi}_1(\xi)\|_{L^2}}\vspace{2mm}\\ & & {\displaystyle
+\int^t_0\|(1+|\xi|)^{\frac{s}{2}}\cos((t-\tau)|\xi|)\hat{f}(\tau,
\xi)\|_{L^2}d\tau}\vspace{2mm}\\ & \leq & {\displaystyle
C\left(\|\partial_x
\phi_0\|_{H^s}+\|\phi_1\|_{H^s}+\int^t_0\|f(\tau,
\cdot)\|_{H^s}\right).}\end{array}\end{equation} Similarly, we have
\begin{equation}
\|\partial_x\phi(t, \cdot)\|_{H^s}\leq C\left(\|\partial_x
\phi_0\|_{H^s}+\|\phi_1\|_{H^s}+\int^t_0\|f(\tau,
\cdot)\|_{H^s}\right).
\end{equation}
Thus, (2.10) comes from (2.15) and (2.16) immediately. This proves
Lemma 2.2. $\quad\quad\quad\blacksquare$

\begin{Lemma}
Let $\phi$ be a solution of the Cauchy problem
\begin{equation}\left\{\begin{array}{l}
{\displaystyle \phi_{tt}-\triangle \phi= \sum^n_{j=0}a_j\partial_jf_j,}\vspace{2mm}\\
t=0: \;\; \phi=0, \quad \phi_t=0\end{array}\right.
\end{equation}
Then
\begin{equation}
\|\phi(t,\cdot)\|_{L^2}\leq
C\left(\sum^n_{j=0}\int^t_0\|f_j(\tau,\cdot)\|_{L^2}d\tau+\|f_0(0,
\cdot)\|_{L^2}\right).
\end{equation}
In particular,  for $n\geq 2$ it holds that
\begin{equation}\begin{array}{lll}
|\phi(t, x)| & \leq & {\displaystyle
C(1+t)^{-\frac{n-1}{2}}\left\{\int^t_0(1+\tau)^{\frac{n-1}{2}}\sum^{n}_{j=0}\|f_j(\tau,
\cdot)\|_{L^\infty}d\tau \right.}\vspace{2mm}\\
& & +{\displaystyle
\left.\int^t_0(1+\tau)^{-\frac{n+1}{2}}\sum^{n}_{j=0}\sum_{|I|\leq
n+1}\|Z^If_j(\tau, \cdot)\|_{L^1}d\tau\right\}.}\end{array}
\end{equation}
\end{Lemma}
\noindent {\bf Proof.} Taking the Fourier transformation on the
variable $x$ in (2.17) yields
\begin{equation}\left\{\begin{array}{l}
{\displaystyle \hat{\phi}_{tt}+|\xi|^2\hat{\phi}=\sum^n_{j=1}\sqrt{-1}a_j\xi_j\hat{f}_j+a_0\partial_t\hat{f}_0,}\vspace{2mm}\\
t=0:\;\; \hat{\phi}=0, \quad \hat{\phi}_t=0.\end{array}\right.
\end{equation}
Solving the initial value problem (2.20) gives
\begin{equation}
\hat{\phi}(t,
\xi)=\sum^n_{j=1}a_j\int^t_0\frac{\sin((t-\tau)|\xi|)}{|\xi|}\sqrt{-1}\xi_j\hat{f}_jd\tau+a_0\int^t_0\frac{\sin((t-\tau)|\xi|)}{|\xi|}\partial_t\hat{f}_0d\tau.
\end{equation}
By Minkowski inequality, we have
\begin{equation}
\left\|\sum^n_{j=1}a_j\int^t_0\frac{\sin((t-\tau)|\xi|)}{|\xi|}\sqrt{-1}\xi_j\hat{f}_jd\tau\right\|_{L^2}
\leq C\sum^n_{j=1}\|f_j(\tau,\cdot)\|_{L^2}d\tau.
\end{equation}
Using the integration by parts, we obtain
$$\begin{array}{lll}
{\displaystyle
a_0\int^t_0\frac{\sin((t-\tau)|\xi|)}{|\xi|}\partial_t\hat{f}_0d\tau}
& = & {\displaystyle
a_0\int^t_0\frac{\sin((t-\tau)|\xi|)}{|\xi|}d\hat{f}_0}\vspace{2mm}\\
& = & {\displaystyle  -a_0\frac{\sin(t|\xi|)}{|\xi|}\hat{f}_0(0,
\xi)+a_0\int^t_0\cos((t-\tau)|\xi|)\hat{f}_0(\tau,
\xi)d\tau.}\end{array}$$ It follows from the Minkowski inequality
that
\begin{equation}\begin{array}{l}
{\displaystyle
\left\|a_0\int^t_0\frac{\sin((t-\tau)|\xi|)}{|\xi|}\partial_t\hat{f}_0d\tau\right\|_{L^2}}\vspace{2mm}\\
\qquad\leq  {\displaystyle
|a_0|\left\|\frac{\sin(t|\xi|)}{|\xi|}\hat{f}_0(0,
\xi)\right\|_{L^2}+|a_0|\left\|\int^t_0\cos((t-\tau)|\xi|)\hat{f}_0(\tau,
\xi)d\tau\right\|_{L^2}}\vspace{2mm}\\
\qquad \leq  {\displaystyle C\|f(0,
\cdot)\|_{\dot{H}^{-1}}+C\int^t_0\|f_0(\tau,
\cdot)\|_{L^2}d\tau.}\end{array}
\end{equation}
Noting the definition of $\dot{H}^{-1}$ and using H\"{o}lder
inequality, we have
\begin{equation}
\|f(0, \cdot)\|_{\dot{H}^{-1}}=\sup_{v\in H^1, v\neq
0}\frac{\int_{\mathbb{R}^n}f(0, \xi)v(\xi)d\xi}{\|v\|_{H^1}}\leq
\|f(0, \cdot)\|_{L^2}.
\end{equation}
Combining (2.23) and (2.24) yields
\begin{equation}
\left\|a_0\int^t_0\frac{\sin((t-\tau)|\xi|)}{|\xi|}\partial_t\hat{f}_0d\tau\right\|_{L^2}\leq
C\|f(0, \cdot)\|_{L^2}+C\int^t_0\|f_0(\tau, \cdot)\|_{L^2}d\tau.
\end{equation}
Thus, we obtain (2.18) immediately from (2.21), (2.22), (2.25) and
Minkowski inequality.

The proof of (2.19) can be found in Li and Zhou \cite{lz1}, here we
omit it. Thus the proof of Lemma 2.3 is completed.
$\quad\quad\quad\blacksquare$

The following lemma comes from Klainerman \cite{kl1}.
\begin{Lemma} Suppose that $\phi$ is $C^2$ smooth and satisfies
$$\square \phi+\sum^n_{j,k=0}\gamma^{jk}(t, x)\partial_j\partial_k \phi=F \quad  (0\leq t\leq T),$$
and suppose furthermore that
$$\phi\longrightarrow 0\quad{\rm as}\;\;|x|\rightarrow \infty.$$
If
$$|\gamma|=\sum^n_{j,k=0}|\gamma^{jk}|\leq \frac{1}{2} \quad  (0\leq t\leq T),$$
then, for any given $t\in [0,T]$, it holds that
\begin{equation}
\|\partial \phi(t, \cdot)\|_{L^2}\leq
2\exp\left\{\int^t_02|\dot{\gamma(\tau)}|d\tau\right\}\|\partial
\phi(0,
\cdot)\|_{L^2}+2\int^t_0\exp\left\{\int^t_s2|\dot{\gamma(\tau)}|d\tau\right\}\|F(s,
\cdot)\|_{L^2}ds,
\end{equation}
where $$|\dot{\gamma}(t)|=\sup|\partial_i\gamma^{jk}(t, \cdot)|.$$
\end{Lemma}

\begin{Lemma}
Suppose that $ G=G(w)$  is a sufficiently smooth function of $
w=(w_1, \cdots, w_m)$ with
\begin{equation}
G(0)=0.
\end{equation}
For any given integer $N \geq 0$, if a vector function $w=w(t, x)$
satisfies
\begin{equation}
\sum_{|I|\leq [\frac{N}{2}]}\|Z^Iw(t, \cdot)\|_{L^\infty}\leq \nu_0,
\quad \forall\; t\in [0, T],
\end{equation}
where $[\cdot]$ stands for the integer part of a real number and $
\nu_0 $ is a positive constant, then it holds that
\begin{equation}
\sum_{|I|\leq N}\|Z^IG(w(t, \cdot))\|_{L^P}\leq
C(\nu_0)\sum_{|I|\leq N}\|Z^Iw(t,\cdot)\|_{L^p}, \quad \forall\;
t\in [0, T],
\end{equation}
provided that all norms appearing on the right-hand side of (2.29)
are bounded, where $ C(\nu_0)$ is a positive constant depending on $
\nu_0$, and $p$ is a real number with $1\leq p\leq \infty$.
\end{Lemma}

The proof of Lemma 2.5 can be found in Li and Chen \cite{lc}.

\begin{Lemma}
Assume that $I=(I_1,\cdots, I_\sigma)$ and $J=(J_1, \cdots,
J_\sigma)$ is a multi-index. If a vector function $\phi=\phi(t, x)$
satisfies
\begin{equation}
\sum_{|J|\leq [\frac{|I|}{2}]}\|Z^J\phi(t, \cdot)\|_{L^\infty}\leq
\nu_0, \quad \forall\; t\in [0, T],
\end{equation}
then it holds that
\begin{equation}\begin{array}{lll}
{\displaystyle \|Z^I((e^{-\phi}-1)\partial_i\phi)(t, \cdot)\|_{L^2}}
& \leq & {\displaystyle C(\nu_0)\sum_{|I_1|\leq |I|}\sum_{|I_2|\leq
[\frac{|I|-1}{2}]}\|Z^{I_1}\phi(t,
\cdot)\|_{L^2}\|Z^{I_2}\partial_i\phi(t,
\cdot)\|_{L^\infty}+}\vspace{2mm}\\
& & {\displaystyle C(\nu_0)\sum_{|I_2|\leq |I|}\sum_{|I_1|\leq
[\frac{|I|}{2}]}\|Z^{I_1}\phi(t,
\cdot)\|_{L^\infty}\|Z^{I_2}\partial_i\phi(t,
\cdot)\|_{L^2}.}\end{array}\end{equation} provided that all norms
appearing on the right-hand side of (2.31) are bounded.
\end{Lemma}

\noindent {\bf Proof.} When $|I|=0,$ by Lemma 2.5 we have
\begin{equation}\begin{array}{lll}\|(e^{-\phi}-1)\partial_i\phi(t, \cdot)\|_{L^2}&
\leq & \|(e^{-\phi}-1)(t,
\cdot)\|_{L^\infty}\|\partial_i\phi(t, \cdot)\|_{L^2}\vspace{2mm}\\
& \leq & C(\nu_0)\|\phi(t, \cdot)\|_{L^\infty}\|\partial_i\phi(t,
\cdot)\|_{L^2}.\end{array}
\end{equation}
For $|I|\geq 1$, it follows from Minkowski inequality and Lemma 2.5
that
\begin{equation}\begin{array}{lll}{\displaystyle \|Z^I((e^{-\phi}-1)\partial_i\phi)(t, \cdot)\|_{L^2}} & \leq &
{\displaystyle C\sum_{|I_1|+|I_2|\leq |I|, |I_1|>|I_2|}
\|Z^{I_1}(e^{-\phi}-1)(t, \cdot)\|_{L^2}\|Z^{I_2}\partial_i\phi(t,
\cdot)\|_{L^\infty}+}\vspace{2mm}\\
& & {\displaystyle C\sum_{|I_1|+|I_2|\leq |I|, |I_1|\leq|I_2|}
\|Z^{I_1}(e^{-\phi}-1)(t,
\cdot)\|_{L^\infty}\|Z^{I_2}\partial_i\phi(t, \cdot)\|_{L^2 }}\vspace{2mm}\\
&\leq & {\displaystyle C(\nu_0)\sum_{|I_1|\leq |I|}\sum_{|I_2|\leq
[\frac{|I|-1}{2}]}\|Z^{I_1}\phi(t,
\cdot)\|_{L^2}\|Z^{I_2}\partial_i\phi(t, \cdot)\|_{L^\infty}+}\vspace{2mm}\\
& & {\displaystyle C(\nu_0)\sum_{|I_2|\leq |I|}\sum_{|I_1|\leq
[\frac{|I|}{2}]}\|Z^{I_1}\phi(t,
\cdot)\|_{L^\infty}\|Z^{I_2}\partial_i\phi(t,
\cdot)\|_{L^2}.}\end{array}\end{equation}
 (2.31) follows from (2.32)
and (2.33) immediately. Thus the proof of Lemma 2.6 is completed.
$\quad\quad\quad\blacksquare$

\begin{Lemma} Suppose that $\phi_0(x), \phi_1(x)\in C^\infty(\mathbb{R}^2)$
and suppose furthermore that there exist two positive constants
$A\in \mathbb{R}^+$ and $k\in \mathbb{R}^+$ such that
$$|\phi_0(x)|\leq \frac{A}{(1+|x|)^k},\quad |\phi_1(x)|\leq
\frac{A}{(1+|x|)^{k+1}}\quad (k>1).\eqno{(H)}$$ If $\phi=\phi(t, x)$
is a solution of the following Cauchy problem
\begin{equation}\left\{\begin{array}{l}
\phi_{tt}-\triangle \phi=0,\vspace{2mm}\\
t=0:\;\; \phi=\phi_0(x), \quad \phi_t=\phi_1(x).\end{array}\right.
\end{equation}
Then it holds that
\begin{equation}|\phi(t, x)|\leq
 \left\{\begin{array}{l}{\displaystyle
\frac{CA}{\sqrt{1+t+|x|}(1+|t-|x||)^{k-\frac{1}{2}}}\quad (|x|\geq
t),}\vspace{3mm}\\
{\displaystyle \frac{CA}{\sqrt{1+t+|x|}\sqrt{1+|t-|x||}}\quad
(|x|\leq t).}
\end{array}\right.
\end{equation}
\end{Lemma}

\begin{Remark} Here we would like to mention that, if the condition (H) is replaced by
$$|\phi_0(x)|\leq \frac{A}{(1+|x|)^{k+1}},\quad |\phi_1(x)|\leq
\frac{A}{(1+|x|)^{k+1}}\quad (k>1).\eqno{(H')}$$  Tsuyata \cite{ts}
has showed that the solution of the Cauchy problem (2.34) satisfies
the following decay estimate
$$|\phi(t, x)|\leq \frac{CA}{\sqrt{1+t+|x|}\sqrt{1+|t-|x||}}.$$
Obviously, Lemma 2.7 improve the Tsuyata's result given in
\cite{ts}.
\end{Remark}

\noindent {\bf Proof of Lemma 2.7.} It is easy to see that the
solution of (2.34) reads
\begin{equation}
\phi(t, x)=\frac{1}{2\pi t^2}\int_{|x-y|\leq
t}\frac{t\phi_0(y)+t^2\phi_1(y)+t\nabla
\phi_0(y)\cdot(y-x)}{(t^2-|y-x|^2)^{\frac{1}{2}}}dy.
\end{equation}

We first estimate ${\displaystyle |\frac{1}{2\pi t}\int_{|x-y|\leq
t}\frac{\phi_0(y)}{(t^2-|y-x|^2)^{\frac{1}{2}}}dy|}$.

Introduce
$$x=(|x|\cos\theta,\; |x|\sin\theta),\quad y=(r\cos(\theta+\psi),\; r\sin(\theta+\psi))$$ and let $\chi$ be the
characteristic function of positive numbers. Then
\begin{equation}\begin{array}{l}{\displaystyle
\left|\frac{1}{2\pi t}\int_{|x-y|\leq
t}\frac{\phi_0(y)}{(t^2-|y-x|^2)^{\frac{1}{2}}}dy\right|}\vspace{3mm}\\
{\displaystyle \qquad \leq \frac{A}{2\pi t}\int_{|x-y|\leq
t}\frac{1}{\sqrt{t^2-|y-x|^2}(1+|y|)^k}dy}\vspace{3mm}\\
\qquad {\displaystyle \leq \frac{A}{2\pi t}
\left(\int^{t+|x|}_{|t-|x||}\frac{r}{(1+r)^k}\int^\varphi_{-\varphi}\frac{1}{\sqrt{t^2-|x|^2-r^2+2r|x|\cos\psi}}d\psi
dr+\right.}\vspace{3mm}\\
 \qquad \quad {\displaystyle\left.
\chi(t-|x|)\int^{t-|x|}_0\frac{r}{(1+r)^k}\int^\pi_{-\pi}\frac{1}{\sqrt{t^2-|x|^2-r^2+2r|x|\cos\psi}}d\psi
dr\right),}\end{array}\end{equation} where
$$\varphi=\arccos\frac{|x|^2+r^2-t^2}{2|x|r}.$$

Let $ h(y) $ be a continuous function on $\mathbb{R}$ and
$y=(r\cos(\theta+\psi),\; r\sin(\theta+\psi))$. Define
$$H(t, |x|, r, \theta, h)=
\left\{\begin{array}{l}{\displaystyle
\int^\varphi_{-\varphi}\frac{h(r,\theta+\psi
)}{\sqrt{t^2-|x|^2-r^2+2|x|r\cos\psi}}d\psi,
\qquad \left|\frac{|x|^2+r^2-t^2}{2|x|r}\right|\leq1},\vspace{3mm}\\
{\displaystyle \int^\pi_{-\pi}\frac{h(r,\theta+\psi
)}{\sqrt{t^2-|x|^2-r^2+2|x|r\cos\psi}}d\psi, \qquad
\left|\frac{|x|^2+r^2-t^2}{2|x|r}\right|\geq1}
\end{array}\right.$$
and
$$H(t, |x|, r)=H(t, |x|, r, \theta, 1),$$
where, as before, $\varphi$ is given by
$$\varphi=\arccos\frac{|x|^2+r^2-t^2}{2|x|r}.$$

The following proposition has been proved in Kovalyov \cite{ko}.
\begin{Proposition}
(I)\;\;  If
$$t\geq |x|+r\quad {\rm and} \quad
\left|\frac{|x|^2+r^2-t^2}{2|x|r}\right|\geq1,$$ then $H(t, |x|, r)$
satisfies
\begin{equation}
H(t, |x|, r)\leq
C\frac{\ln\left\{2+\frac{r|x|}{t^2-(r+|x|)^2}\right\}}{\sqrt{t^2-|x|^2-r^2}}\leq
\frac{C}{t^2-(r+|x|)^2},
\end{equation}
here and hereafter $C$ stands for some constants.

(II)\;\; If
$$t\leq |x|+r \quad {\rm and}\quad
\left|\frac{|x|^2+r^2-t^2}{2|x|r}\right|\leq1,$$ then
\begin{equation}
H(t, |x|, r)\leq
\frac{C}{\sqrt{r|x|}}\ln\left\{2+\frac{r|x|\chi(t-|x|)}{(r+|x|)^2-t^2}\right\},
\end{equation}
where $ \chi$ is the characteristic function of positive numbers.
\end{Proposition}

We now continue to estimate (2.37).

To do so, we distinguish the following two cases: $|x|\geq t$ and
$|x|\leq t$.

\vskip 3mm

{\bf Case I: $|x|\geq t$}

It follows from (2.39) that
\begin{equation}
\left|\frac{1}{2\pi t}\int_{|x-y|\leq
t}\frac{\phi_0(y)}{(t^2-|y-x|^2)^{\frac{1}{2}}}dy\right|\leq
\frac{CA}{t\sqrt{|x|}}\int^{t+|x|}_{|x|-t}\frac{1}{(1+r)^{k-\frac{1}{2}}}dr.
\end{equation}
In the present situation, we distinguish the following cases:
$t\geq1$ and $0<t<1$.

\vskip 2mm

{\bf Case I-A: $t\geq1$}

In this case, according to $k$, we distinguish the following three
cases:

 \vskip 2mm

{\bf Case I-A-1: $k> \frac{3}{2}$}

In the present situation, it holds that
$$\frac{CA}{t\sqrt{|x|}}\int^{t+|x|}_{|x|-t}\frac{1}{(1+r)^{k-\frac{1}{2}}}dr
=\frac{CA}{t\sqrt{|x|}(1+|x|-t)^{k-\frac{3}{2}}}\left[1-\left(\frac{1+|x|-t}{1+|x|+t}\right)^{k-\frac{3}{2}}\right].$$
Noting
$$1-s^{k-\frac{3}{2}}\leq C(1-s), \quad \forall\; s\in [0,1]$$
and
$$1-\frac{1+|x|-t}{1+|x|+t}=\frac{2t}{1+|x|+t},$$
we have
\begin{equation}\begin{array}{lll}{\displaystyle \left|\frac{1}{2\pi
t}\int_{|x-y|\leq
t}\frac{\phi_0(y)}{(t^2-|y-x|^2)^{\frac{1}{2}}}dy\right|} & \leq &
{\displaystyle
\frac{CA}{\sqrt{|x|}(1+|x|-t)^{k-\frac{3}{2}}(1+|x|+t)}},\vspace{3mm}\\
& \leq & {\displaystyle \frac{
CA}{\sqrt{|x|+t}(1+|x|-t)^{k-\frac{1}{2}}}}.\end{array}
\end{equation}

\vskip 2mm

{\bf Case I-A-2: $k=\frac{3}{2}$}

It follows from (2.40) that
\begin{equation}\begin{array}{lll}{\displaystyle \left|\frac{1}{2\pi t}\int_{|x-y|\leq
t}\frac{\phi_0(y)}{(t^2-|y-x|^2)^{\frac{1}{2}}}dy\right|}
 & \leq & {\displaystyle
 \frac{CA}{t\sqrt{|x|}}\int^{t+|x|}_{|x|-t}\frac{1}{(1+r)}dr
 =\frac{CA}{t\sqrt{|x|}}\ln\left\{1+\frac{2t}{1+|x|-t}\right\}}\vspace{3mm}\\
& \leq & {\displaystyle \frac{ CA}{\sqrt{|x|}(1+|x|-t)}\leq
 \frac{ CA}{\sqrt{|x|+t}(1+|x|-t)}.}\end{array}\end{equation}

\vskip 2mm

{\bf Case I-A-3: $1< k <\frac{3}{2}$}

In the present situation, it follows from (2.40) that
$$\begin{array}{lll}{\displaystyle \left|\frac{1}{2\pi t}\int_{|x-y|\leq
t}\frac{\phi_0(y)}{(t^2-|y-x|^2)^{\frac{1}{2}}}dy\right|} & \leq &
{\displaystyle
\frac{CA}{t\sqrt{|x|}}\left[(1+t+|x|)^{\frac{3}{2}-k}-(1+|x|-t)^{\frac{3}{2}-k}\right]
}\vspace{3mm}\\
& = & {\displaystyle
\frac{CA}{t\sqrt{|x|}(1+|x|-t)^{k-\frac{3}{2}}}\left[\left(\frac{1+t+|x|}{1+|x|-t}\right)^{\frac{3}{2}-k}-1\right].}
\end{array}$$
Noting the fact that $1< k<\frac{3}{2}$, we have
$$\left(\frac{1+t+|x|}{1+|x|-t}\right)^{\frac{3}{2}-k}-1\leq \frac{C t}{1+|x|-t}.$$
Hence,
\begin{equation}
\left|\frac{1}{2\pi t}\int_{|x-y|\leq
t}\frac{\phi_0(y)}{(t^2-|y-x|^2)^{\frac{1}{2}}}dy\right|\leq \frac{
CA}{\sqrt{|x|+t}(1+|x|-t)^{k-\frac{1}{2}}}.
\end{equation}

Summarizing the above argument, for the case that $|x|\geq t$ and
$t\geq 1$, we obtain from (2.41)-(2.43) that
\begin{equation}
\left|\frac{1}{2\pi t}\int_{|x-y|\leq
t}\frac{\phi_0(y)}{(t^2-|y-x|^2)^{\frac{1}{2}}}dy\right|\leq\frac{
CA}{\sqrt{1+|x|+t}(1+|x|-t)^{k-\frac{1}{2}}} \quad (k>1).
\end{equation}

\vskip 2mm

{\bf Case I-B: $|x|\geq t$ and $0<t<1$}

We next consider the case that $|x|\geq t$ and $0<t <1$. In this
case, we distinguish the following two cases.

\vskip 2mm

{\bf Case I-B-1: $|t-|x||\leq 1$}

Introducing the variable $r=|x-y|$, we have
\begin{equation}\begin{array}{lll}{\displaystyle \left|\frac{1}{2\pi t}\int_{|x-y|\leq
t}\frac{\phi_0(y)}{(t^2-|y-x|^2)^{\frac{1}{2}}}dy\right|} & \leq &
{\displaystyle \frac{1}{2\pi t}\int_{|x-y|\leq
t}\frac{CA}{(t^2-|y-x|^2)^{\frac{1}{2}}(1+|y|)^k}dy}\vspace{3mm}\\
& \leq & {\displaystyle \frac{CA}{\pi
t}\int^t_0\frac{r}{\sqrt{t^2-r^2}}dr \leq
\frac{CA}{\sqrt{1+t+|x|}(1+|x|-t)^{k-\frac{1}{2}}}.}\end{array}
\end{equation}

\vskip 2mm

{\bf Case I-B-2: $|t-|x||>1$}

Noting the fact that $|x|\geq t$ and $0<t<1$, we observe
$$|x|>t+1.$$
Thus, by the case (I-A), we have
\begin{equation}\left|\frac{1}{2\pi t}\int_{|x-y|\leq
t}\frac{\phi_0(y)}{(t^2-|y-x|^2)^{\frac{1}{2}}}dy\right| \leq
\frac{CA}{\sqrt{1+t+|x|}(1+|x|-t)^{k-\frac{1}{2}}}.
\end{equation}
Therefore, combining (2.44)-(2.46) gives
\begin{equation}\left|\frac{1}{2\pi t}\int_{|x-y|\leq
t}\frac{\phi_0(y)}{(t^2-|y-x|^2)^{\frac{1}{2}}}dy\right|\leq
\frac{CA}{\sqrt{1+t+|x|}(1+|x|-t)^{k-\frac{1}{2}}},
\end{equation}
provided that $|x|\geq t$.

\vskip 3mm

{\bf Case II: $|x|\leq t$}

We now consider the case that $|x|\leq t$.

It follows from (2.37) that
\begin{equation}
\left|\frac{1}{2\pi t}\int_{|x-y|\leq
t}\frac{\phi_0(y)}{(t^2-|y-x|^2)^{\frac{1}{2}}}dy\right|\leq I+II,
\end{equation}
where
$$I=\frac{A}{2\pi t}\int^{t+|x|}_{t-|x|}\frac{H(t, |x|, r)r}{(1+r)^k}, \quad
II=\frac{A}{2\pi t}\int^{t-|x|}_{0}\frac{H(t, |x|, r)r}{(1+r)^k}.$$

We next estimate $I$ and $II$ by distinguishing the follows cases.

\vskip 2mm

{\bf Case II-A: $t+|x|\geq 1$}

It follows from (2.39) that
\begin{equation}
I\leq\frac{CA}{t\sqrt{|x|}}\int^{t+|x|}_{t-|x|}\ln\left\{2+\frac{|x|}{|x|+r-t}\right\}\frac{1}{(1+r)^{k-\frac{1}{2}}}dr.
\end{equation}
Introducing the variable $\xi=|x|+r-t$, we have
\begin{equation}\begin{array}{lll}{\displaystyle I} & \leq &{\displaystyle \frac{CA}{t\sqrt{|x|}}\int^{2|x|}_{0}
\ln\left\{2+\frac{|x|}{\xi}\right\}\frac{1}{(1+\xi+t-|x|)^{k-\frac{1}{2}}}d\xi}\vspace{3mm}\\
& \leq & {\displaystyle
\frac{CA}{t\sqrt{|x|}(1+t-|x|)^{k-\frac{1}{2}}}\int^{2|x|}_{0}\ln\left\{\frac{3|x|}{\xi}\right\}d\xi}\vspace{3mm}\\
& = & {\displaystyle
\frac{CA}{t\sqrt{|x|}(1+t-|x|)^{k-\frac{1}{2}}}[2|x|\ln(3|x|)-2|x|\ln(2|x|)+2|x|]}\vspace{3mm}\\
& \leq & {\displaystyle
\frac{CA}{\sqrt{t}(1+t-|x|)^{k-\frac{1}{2}}}\leq
\frac{CA}{\sqrt{t+|x|}(1+t-|x|)^{k-\frac{1}{2}}}}\vspace{3mm}\\
& \leq & {\displaystyle
\frac{CA}{\sqrt{1+t+|x|}(1+t-|x|)^{k-\frac{1}{2}}}.}\end{array}
\end{equation}

We now estimate $II$.

By (2.38), we have
$$\begin{array}{lll}{\displaystyle II} & \leq & {\displaystyle
\frac{CA}{t}\int^{t-|x|}_{0}\frac{1}{\sqrt{t^2-(|x|+r)^2}(1+r)^{k-1}}dr}\vspace{3mm}\\
& \leq & {\displaystyle
\frac{CA}{t\sqrt{t+|x|}}\int^{t-|x|}_{0}\frac{1}{\sqrt{t-|x|-r}(1+r)^{k-1}}dr.}\end{array}$$
Let $$\rho=\sqrt{t-|x|-r}.$$ Then
\begin{equation}\begin{array}{lll}{\displaystyle
II} & \leq & {\displaystyle
\frac{CA}{t\sqrt{t+|x|}}\int^{\sqrt{t-|x|}}_{0}
\frac{1}{(1+t-|x|-\rho^2)^{k-1}}d\rho }\vspace{3mm}\\
& \leq & {\displaystyle
\frac{CA}{t\sqrt{t+|x|}(1+t-|x|)^{\frac{k-1}{2}}}\int^{\sqrt{t-|x|}}_{0}\frac{1}{(\sqrt{1+t-|x|}-\rho)^{k-1}}d\rho.}
\end{array}\end{equation}
In order to estimate $II$, we distinguish the following three cases.

\vskip 2mm

{\bf Case II-A-1: $k>2 $}

In the present situation, it follows from (2.51) that
$$\begin{array}{lll}{\displaystyle II} &
\leq &{\displaystyle
\frac{CA}{t\sqrt{t+|x|}(1+t-|x|)^{\frac{k-1}{2}}}\left\{\frac{1}{(\sqrt{1+t-|x|}-\sqrt{t-|x|})^{k-2}}-
(\sqrt{1+t-|x|})^{k-2}\right\}}\vspace{3mm}\\
& \leq & {\displaystyle
\frac{CA}{t\sqrt{t+|x|}(1+t-|x|)^{\frac{k-1}{2}}}(\sqrt{1+t-|x|}+\sqrt{t-|x|})^{k-2}}\vspace{3mm}\\
& \leq & {\displaystyle \frac{CA}{t\sqrt{t+|x|}\sqrt{1+t-|x|}}\leq
\frac{CA}{\sqrt{1+t+|x|}\sqrt{1+t-|x|}}.}\end{array}$$

\vskip 2mm

{\bf Case II-A-2: $k=2$}

In this case, by (2.51) we have
$$\begin{array}{lll}{\displaystyle II} &
\leq &{\displaystyle
\frac{CA}{t\sqrt{t+|x|}\sqrt{1+t-|x|}}\times\ln\left\{\frac{\sqrt{1+t-|x|}}{\sqrt{1+t-|x|}-\sqrt{t-|x|}}\right\}}\vspace{3mm}\\
& \leq & {\displaystyle
\frac{CA}{t\sqrt{t+|x|}\sqrt{1+t-|x|}}\times \frac{\sqrt{1+t-|x|}}{\sqrt{1+t-|x|}-\sqrt{t-|x|}}}\vspace{3mm}\\
& \leq & {\displaystyle
\frac{CA}{\sqrt{1+t+|x|}\sqrt{1+t-|x|}}.}\end{array}$$

\vskip 2mm

{\bf Case II-A-3: $1<k<2$}

In this situation, we obtain from (2.51) that
$$\begin{array}{lll}{\displaystyle II} &
\leq &{\displaystyle
\frac{CA}{t\sqrt{t+|x|}(1+t-|x|)^{\frac{k-1}{2}}}\left\{(1+t-|x|)^{\frac{-k+2}{2}}-
(\sqrt{1+t-|x|}-\sqrt{t-|x|})^{-k+2}\right\}}\vspace{3mm}\\
& \leq & {\displaystyle
\frac{CA}{\sqrt{1+t+|x|}\sqrt{1+t-|x|}}.}\end{array}$$

Summarizing the above argument gives
\begin{equation}
II\leq \frac{CA}{\sqrt{1+t+|x|}\sqrt{1+t-|x|}},
\end{equation}
provided that $t+|x|\geq 1$.

\vskip 2mm

{\bf Case II-B: $0<t+|x|< 1$}

As before, introducing the variable $r=|x-y|$, we have
\begin{equation}\begin{array}{lll}{\displaystyle
\left|\frac{1}{2\pi t}\int_{|x-y|\leq
t}\frac{\phi_0(y)}{(t^2-|y-x|^2)^{\frac{1}{2}}}dy\right|} & \leq &
{\displaystyle \frac{1}{2\pi t}\int_{|x-y|\leq
t}\frac{CA}{(t^2-|y-x|^2)^{\frac{1}{2}}(1+|y|)^k}dy}\vspace{3mm}\\
& \leq & {\displaystyle \frac{CA}{\pi
t}\int^t_0\frac{r}{\sqrt{t^2-r^2}}dr \leq
\frac{CA}{\sqrt{1+t+|x|}\sqrt{1+t-|x|}}.}\end{array}
\end{equation}
Combining (2.50) and (2.52)-(2.53) leads to
\begin{equation}
\left|\frac{1}{2\pi t}\int_{|x-y|\leq
t}\frac{\phi_0(y)}{(t^2-|y-x|^2)^{\frac{1}{2}}}dy\right|\leq
\frac{CA}{\sqrt{1+t+|x|}\sqrt{1+t-|x|}}\quad(|x|\leq t).
\end{equation}
(2.47) and (2.54) imply
\begin{equation}
\left|\frac{1}{2\pi t}\int_{|x-y|\leq
t}\frac{\phi_0(y)}{(t^2-|y-x|^2)^{\frac{1}{2}}}dy\right|\leq
\left\{\begin{array}{l}{\displaystyle
\frac{CA}{\sqrt{1+t+|x|}(1+|t-|x||)^{k-\frac{1}{2}}}\quad (|x|\geq
t),}\vspace{3mm}\\
{\displaystyle \frac{CA}{\sqrt{1+t+|x|}\sqrt{1+|t-|x||}}\quad
(|x|\leq t).}
\end{array}\right.
\end{equation}
By Tsutaya \cite {ts}, we have
\begin{equation}
\left|\frac{1}{2\pi }\int_{|x-y|\leq
t}\frac{\phi_1(y)}{(t^2-|y-x|^2)^{\frac{1}{2}}}dy\right|\leq
\left\{\begin{array}{l}{\displaystyle
\frac{CA}{\sqrt{1+t+|x|}(1+|t-|x||)^{k-\frac{1}{2}}}\quad (|x|\geq
t),}\vspace{3mm}\\
{\displaystyle \frac{CA}{\sqrt{1+t+|x|}\sqrt{1+|t-|x||}}\quad
(|x|\leq t)}
\end{array}\right.
\end{equation}
and
\begin{equation}
\left|\frac{1}{2\pi t}\int_{|x-y|\leq t}\frac{\nabla
\phi_0(y)\cdot(y-x)}{(t^2-|y-x|^2)^{\frac{1}{2}}}dy\right|\leq
\left\{\begin{array}{l}{\displaystyle
\frac{CA}{\sqrt{1+t+|x|}(1+|t-|x||)^{k-\frac{1}{2}}}}\quad (|x|\geq
t),\vspace{3mm}\\
{\displaystyle \frac{CA}{\sqrt{1+t+|x|}\sqrt{1+|t-|x|}}\quad
(|x|\leq t).}
\end{array}\right.
\end{equation}
(2.35) follows from (2.55)-(2.57) and (2.36) immediately. Thus, the
proof of Lemma 2.7 is completed. $\quad\quad\quad\blacksquare$

\begin{Lemma}
Suppose that $ \phi $ is a solution to the Cauchy problem
$$\phi_{tt}-\triangle \phi=g$$
with zero initial data. Then
\begin{equation}
|\phi(t, x)|\leq C(1+t+|x|)^{-\frac{n-1}{2}}\sum_{|I|\leq
n-1}\int^t_0\|(Z^Ig)(\tau,
\cdot)/(1+\tau+|\cdot|)^{\frac{n-1}{2}}\|_{L^1}d\tau.
\end{equation}
In particular, for $n=2$ and $p\in (1,2]$ it holds that
\begin{equation}
\|\phi(t, \cdot)\|_{L^p(\mathbb{R}^2)}\leq
C(1+t)^{\frac{2}{p}-1}\int^t_0\|g(\tau,
\cdot)\|_{L^1(\mathbb{R}^2)}d\tau.
\end{equation}
\end{Lemma}
\noindent {\bf Proof.} The inequality (2.58) comes from H\"{o}mander
\cite{ho1} or Klainerman \cite{kl2} directly, while the proof of
(2.59) has been proved by Li and Zhou \cite{lz}.
$\quad\quad\quad\blacksquare$

\begin{Lemma}
Suppose that $n=2$, and suppose furthermore that $\phi=\phi(t, x)$
is a solution of the wave equation
\begin{equation}
\phi_{tt}-\Delta \phi=|g_1g_2(t, x)|
\end{equation}
with zero initial data. Then it holds that
\begin{equation}
\|\phi(t, \cdot)\|_{L^2(\mathbb{R}^2)}\leq
C(1+t)^{\frac{1}{4}}\left\{\sum_{|I|\leq
1}\int^t_0(1+\tau)^{-\frac{1}{2}}\|\Gamma g_1(\tau,
\cdot)\|^2_{L^2(\mathbb{R}^2)}d\tau\right\}^{\frac{1}{2}}
\left\{\int^t_0\|g_2(\tau,
\cdot)\|^2_{L^2(\mathbb{R}^2)}d\tau\right\}^{\frac{1}{2}}
\end{equation}
and
\begin{equation}
(1+t)^{\frac{1}{2}}\|\phi(t, \cdot)\|_{L^\infty(\mathbb{R}^2)}\leq
C\left\{\int^t_0\sum_{|I|\leq 1}\|\Gamma^Ig_1(\tau, \cdot)
\|^2_{L^2}\frac{d\tau}{\sqrt{1+\tau}}\right\}^{\frac{1}{2}}\left\{\int^t_0\sum_{|I|\leq
1}\|\Gamma^Ig_2(\tau, \cdot)
\|^2_{L^2}\frac{d\tau}{\sqrt{1+\tau}}\right\}^{\frac{1}{2}}.
\end{equation}
\end{Lemma}

\noindent {\bf Proof.} The proof of (2.61) can be found in
\cite{lz}. In what follows, we prove (2.62).

Let $E$ be the forward fundamental solution of the wave operator. By
the positivity of $E$ and the H\"{o}lder inequality, we have
\begin{equation}
\phi(t, x)\leq E*|g_1g_2(t, x)|\leq \left(E*g^2_1(t,
x)\right)^{\frac{1}{2}}\left(E*g^2_2(t, x)\right)^{\frac{1}{2}}.
\end{equation}
It follows from Lemma 2.8 and H\"{o}lder inequality that
\begin{equation}
E*g^2_1(t, x)\leq C(1+t)^{-\frac{1}{2}}\sum_{|I|\leq
1}\int^t_0\|Z^Ig_1(\tau,
\cdot)\|^2_{L^2}\frac{d\tau}{\sqrt{1+\tau}}.
\end{equation}
Similarly,
\begin{equation}
E*g^2_2(t, x)\leq C(1+t)^{-\frac{1}{2}}\sum_{|I|\leq
1}\int^t_0\|Z^Ig_2(\tau,
\cdot)\|^2_{L^2}\frac{d\tau}{\sqrt{1+\tau}}.
\end{equation}
(2.62) comes from (2.63)-(2.65) immediately. This proves Lemma 2.9.
$\quad\quad\quad\blacksquare$

The following lemma can be found in Klainerman \cite{kl3}.
\begin{Lemma}
Assume that $p\in [1,\infty)$ and $N$ is an integer satisfying $N
>\frac{n}{p}$. Then it holds that
\begin{equation}
|\phi(t, x)|\leq
C(1+t+|x|)^{-\frac{n-1}{p}}(1+|t-|x||)^{-\frac{1}{p}}\sum_{|I|\leq
N}\|Z^I\phi(t,\cdot)\|_{L^p},
\end{equation}
provided that all norms appearing on the right-hand side of (2.66)
are bounded.
\end{Lemma}

\section{Lower bound of life-span}

This section is devoted to the proof of Theorem 1.1. In order to
prove Theorem 1.1, it suffices to show the following theorem.

\begin{Theorem} Suppose that $u_0(x),\; u_1(x)\in C^\infty(\mathbb{R}^2)$ and
satisfy that there exist two positive constants $A\in \mathbb{R}^+$
and $k\in \mathbb{R}^{+}$ such that $$|u_0(x)|\leq
\frac{A}{(1+|x|)^k},\quad |u_1(x)|\leq \frac{A}{(1+|x|)^{k+1}}\quad
(k>1).$$ Then there exist two positive constants $\delta $ and
$\varepsilon_0$ such that for any fixed $ \varepsilon\in
(0,\varepsilon_0]$, the Cauchy problem (\ref{1.8})-(\ref{1.9}) has a
unique $C^\infty$ solution on the interval $[0, T_\varepsilon]$,
where $T_\varepsilon$ is given by
\begin{equation}\label{1.11}
T_\varepsilon=\frac{\delta}{\varepsilon^{\frac{4}{3}}}-1.
\end{equation}
\end{Theorem}
\noindent {\bf Proof.} The local existence of classical solutions
has been proved by the method of Picard iteration (see Sogge
\cite{sc} and H\"{o}rmander \cite{ho2}). In what follows, we prove
Theorem 3.1 by the method of continuous induction, or say, the
bootstrap argument.

Let $l_1$ and $l_2$ be two positive integers such that $$l_1-3\geq
l_2\geq \frac{1}{2}[l_1]+1.$$

Introduce
\begin{equation}\left\{\begin{array}{l}
{\displaystyle M_1(t)=\sum_{|I|\leq l_1}\|\partial Z^Iu(t,
\cdot)\|_{L^2(\mathbb{R}^2)},}\vspace{2mm}\\
{\displaystyle
 M_{2}(t)=\sum_{|I|\leq l_1}\| Z^Iu(t, \cdot)\|_{L^2(\mathbb{R}^2)},}\vspace{2mm}\\
{\displaystyle
 N_{1}(t)=\sum_{|I|\leq l_2}\|\partial Z^Iu(t,
\cdot)\|_{L^\infty(\mathbb{R}^2)}},\vspace{2mm}\\
{\displaystyle N_2(t)=\sum_{|I|\leq l_2 }\| Z^Iu(t,
\cdot)\|_{L^\infty(\mathbb{R}^2)}.}\end{array}\right.\end{equation}
By the bootstrap argument, for the time being it is supposed that
there exist some positive constants $M_i,\; N_i\;(i=1, 2)$ and $\mu$
such that
\begin{equation}\left\{\begin{array}{l}
{\displaystyle M_1(t)\leq M_1\varepsilon,}\vspace{2mm}\\
{\displaystyle
 M_2(t)\leq M_2\varepsilon (1+t)^{\frac{1}{4}},}\vspace{2mm}\\
{\displaystyle (1+t)^{\frac{1}{2}}N_1(t)\leq N_1\varepsilon,
}\vspace{2mm}\\
{\displaystyle
 (1+t)^{\frac{1}{2}}N_2(t)\leq N_2\varepsilon,}\end{array}\right.\end{equation}
provided that $\varepsilon,\; \mu $ are suitably small and satisfy
$$\varepsilon (1+t)^{\frac{3}{4}}\leq \mu.\eqno{(3.3a)}$$

According to the bootstrap argument, in what follows we show that,
by choosing  $M_i$ and $N_i\;(i=1,2)$ sufficiently large and
$\varepsilon$ suitably small such that
$$\left\{\begin{array}{l}
{\displaystyle M_1(t)\leq \frac{1}{2}M_1\varepsilon,}\vspace{2mm}\\
{\displaystyle
 M_2(t)\leq \frac{1}{2}M_2\varepsilon (1+t)^{\frac{1}{4}},}\vspace{2mm}\\
{\displaystyle (1+t)^{\frac{1}{2}}N_1(t)\leq
\frac{1}{2}N_1\varepsilon,
}\vspace{2mm}\\
{\displaystyle
 (1+t)^{\frac{1}{2}}N_2(t)\leq \frac{1}{2} N_2\varepsilon,}\end{array}\right.\eqno{(3.3b)}$$
provided that $\varepsilon,\; \mu $ are suitably small and (3.3a)
holds.

We first estimate $M_1(t)$.

The equation (\ref{1.8}) can be rewritten as
\begin{equation}
\square u=(e^{-u}-1)\Delta u- u^2_t.
\end{equation}
It follows from Lemma 2.1 and (3.4) that
\begin{equation}\begin{array}{lll}
 \square Z^Iu & = & {\displaystyle \sum_{|I_1|+|I_2|\leq|I|}\sum_{0\leq i_1, i_2\leq
2}A_{II_1I_2}Z^{I_1}(e^{-u}-1)\partial_{i_1i_2}Z^{I_2}u+}\vspace{3mm}\\
& & {\displaystyle \sum_{|I_1|+|I_2|\leq|I|}\sum_{ 0\leq i_1,
i_2\leq
2}\tilde{A}_{II_1I_2}\partial_{i_1}Z^{I_1}u\partial_{i_2}Z^{I_2}u}.\end{array}
\end{equation}
By Minkowski inequality, (3.2) and Lemma 2.5, for $I$ with $|I|\leq
l_1-1$ we have
\begin{equation}\begin{array}{l}
{\displaystyle \left\|\sum_{|I_1|+|I_2|\leq|I| }\sum_{0\leq i_1,
i_2\leq 2}A_{II_1I_2}Z^{I_1}(e^{-u}-1)\partial_{i_1i_2}Z^{I_2}u(t,
\cdot)
\right\|_{L^2}}\vspace{3mm}\\
\qquad\quad {\displaystyle \leq C\sum_{|I_1|+|I_2|\leq|I|,|I_1|
>|I_2|}\sum_{0\leq i_1, i_2\leq
2}\|Z^{I_1}(e^{-u}-1)\partial_{i_1i_2}Z^{I_2}u(t,
\cdot)\|_{L^2}}+\vspace{3mm}\\
\qquad\quad {\displaystyle \;\;\;\; C\sum_{|I_1|+|I_2|\leq|I|,|I_1|
\leq|I_2|}\sum_{0\leq i_1, i_2\leq
2}\|Z^{I_1}(e^{-u}-1)\partial_{i_1i_2}Z^{I_2}u(t, \cdot)\|_{L^2}}\vspace{3mm}\\
\qquad\quad {\displaystyle \leq C\sum_{|I_1|+|I_2|\leq|I|,|I_1|
>|I_2|}\sum_{0\leq i_1, i_2\leq
2}\|Z^{I_1}(e^{-u}-1)\|_{L^2}\|\partial_{i_1i_2}Z^{I_2}u(t,
\cdot)\|_{L^\infty}}+\vspace{3mm}\\
\qquad\quad {\displaystyle\;\;\;\; C\sum_{|I_1|+|I_2|\leq|I|,|I_1|
\leq|I_2|}\sum_{0\leq i_1, i_2\leq
2}\|Z^{I_1}(e^{-u}-1)\|_{L^\infty}\|\partial_{i_1i_2}Z^{I_2}u(t,
\cdot)\|_{L^2}}\vspace{3mm}\\
\qquad\quad {\displaystyle \leq
C\{M_2(t)N_1(t)+M_1(t)N_2(t)\}.}\end{array}
\end{equation}
provided that $\varepsilon_0>0$ is suitably small.

Again by Minkowski inequality and (3.2), for $I$ with $|I|\leq
l_1-1$ we get
\begin{equation}\begin{array}{l}
{\displaystyle \left\|\sum_{|I_1|+|I_2|\leq|I|}\sum_{ 0\leq i_1,
i_2\leq
2}\tilde{A}_{II_1I_2}\partial_{i_1}Z^{I_1}u\partial_{i_2}Z^{I_2}u\right\|_{L^2}}\vspace{3mm}\\
\qquad\quad {\displaystyle \leq C\sum_{|I_1|+|I_2|\leq |I|,
|I_1|>|I_2|}\sum_{ 0\leq i_1, i_2\leq 2}\|\partial_{i_1}Z^{I_1}u(t,
\cdot)\|_{L^2}\|\partial_{i_2}Z^{I_1}u(t, \cdot)\|_{L^\infty}+}\vspace{3mm}\\
\qquad\quad {\displaystyle\;\;\;\; C\sum_{|I_1|+|I_2|\leq|I|,
|I_1|\leq|I_2|}\sum_{ 0\leq i_1, i_2\leq
2}\|\partial_{i_1}Z^{I_1}u(t,
\cdot)\|_{L^\infty}\|\partial_{i_2}Z^{I_1}u(t, \cdot)\|_{L^2}}\vspace{3mm}\\
\qquad\quad{\displaystyle  \leq CM_1(t)N_1(t).}\end{array}
\end{equation}
On the other hand, noting Lemma 2.2 and using (3.5)-(3.7), for $I$
with $|I|\leq l_1-1$ we have
\begin{equation} \|\partial Z^Iu(t,
\cdot)\|_{L^2}\leq C\|\partial Z^Iu(0,
\cdot)\|_{L^2}+C\int^t_0[M_2(\tau)N_1(\tau)+M_1(\tau)N_2(\tau)+M_1(\tau)N_1(\tau)]d\tau.
\end{equation}

We now estimate $\|\partial Z^Iu(t, \cdot)\|_{L^2}\;(|I|=l_1)$.

By Lemma 2.1 and (3.4),
$$\begin{array}{lll}
{\displaystyle \square Z^Iu} & = & {\displaystyle Z^I\square
u+\sum_{|J|\leq|I|-1}A_{IJ}Z^J\square u}\vspace{3mm}\\
& = & {\displaystyle \sum^2_{i,
j=0}(e^{-u}-1)\partial_i\partial_jZ^Iu+\sum_{|I_1|+|I_2|\leq|I|,
|I_2|\leq|I|-1}\sum_{0\leq i_1, i_2\leq 2}
\bar{A}_{II_1I_2}Z^{I_1}(e^{-u}-1)\partial_{i_1i_2}Z^{I_2}u}\vspace{3mm}\\
&  & +{\displaystyle \sum_{|I_1|+|I_2|\leq|I|, }\sum_{0\leq i_1,
i_2\leq
2}\bar{\bar{A}}_{II_1I_2}\partial_{i_1}Z^{I_1}u\partial_{i_2}Z^{I_2}u.}\end{array}$$
Hence
\begin{equation}\begin{array}{lll}
{\displaystyle \square Z^Iu+\sum^2_{i,
j=0}(1-e^{-u})\partial_i\partial_jZ^Iu} & = & {\displaystyle
\sum_{|I_1|+|I_2|\leq|I|, |I_2|\leq|I|-1}\sum_{0\leq i_1, i_2\leq 2}
\bar{A}_{II_1I_2}Z^{I_1}(e^{-u}-1)\partial_{i_1i_2}Z^{I_2}u}\vspace{3mm}\\
&  & +{\displaystyle \sum_{|I_1|+|I_2|\leq|I|, }\sum_{0\leq i_1,
i_2\leq
2}\bar{\bar{A}}_{II_1I_2}\partial_{i_1}Z^{I_1}u\partial_{i_2}Z^{I_2}u.}\end{array}
\end{equation}
Similar to the proof of (3.6), when $\varepsilon_0>0$ is suitably
small, for $I$ with $|I|=l_1$ it holds that
\begin{equation}
\|\sum_{|I_1|+|I_2|\leq|I|, |I_2|\leq|I|-1}\sum_{0\leq i_1, i_2\leq
2}\bar{A}_{II_1I_2}Z^{I_1}(e^{-u}-1)\partial_{i_1i_2}Z^{I_2}u(t,
\cdot)\|_{L^2} \leq C[M_2(t)N_1(t)+M_1(t)N_2(t)].
\end{equation}
Similar to the proof of (3.7),  for $I$ with $|I|=l_1$ we have
\begin{equation}
\| \sum_{|I_1|+|I_2|\leq|I|, }\sum_{0\leq i_1, i_2\leq
2}\bar{\bar{A}}_{II_1I_2}\partial_{i_1}Z^{I_1}u\partial_{i_2}Z^{I_2}u(t,
\cdot)\|_{L^2}\leq CM_1(t)N_1(t).
\end{equation}

On the other hand, because of (3.3), it holds that
\begin{equation} \sum^2_{i,
j=0}|\gamma^{ij}|\leq CN_2\varepsilon < \frac{1}{2},
\end{equation}
provided that $ \varepsilon$ is suitably small, where
$\gamma^{ij}=(1-e^{-u})$. Moreover, for $|\dot{\gamma}(t)|$ defined
as in Lemma 2.4, we have
\begin{equation}
2\int^t_0|\dot{\gamma}(\tau)|d\tau\leq CN_1\mu\leq \ln2,
\end{equation}
provided that $\varepsilon,\; \mu $ are suitably small and (3.3a)
holds. Thus, noting Lemma 2.4 and using (3.9)-(3.13), for $I$ with
$|I|=l_1$ we have
\begin{equation}
\|\partial Z^Iu(t, \cdot)\|_{L^2}\leq 4\|\partial Z^Iu(0,
\cdot)\|_{L^2}+C\int^t_0(M_2(\tau)N_1(\tau)+M_1(\tau)N_2(\tau)+M_1(\tau)N_1(\tau))d\tau.
\end{equation}
Combining (3.8) and (3.14) yields
\begin{equation}
M_1(t)\leq
K_1\varepsilon+C\int^t_0(M_2(\tau)N_1(\tau)+M_1(\tau)N_2(\tau)+M_1(\tau)N_1(\tau))d\tau.
\end{equation}

We next estimate $M_{2}(t)$.

In the present situation, the equation (\ref{1.8}) can be rewritten
as
\begin{equation}
\square
u=\sum^2_{i=1}\partial_i((e^{-u}-1)\partial_iu)-u^2_t-\sum^2_{i=1}\partial_i(e^{-u}-1)\partial_iu.
\end{equation}
By Lemma 2.1 and (3.16), we obtain
$$\begin{array}{lll}
{\displaystyle \square Z^Iu} & = & {\displaystyle
\sum_{|I_1|+|I_2|\leq |I|}\sum_{0\leq i_1, i_2\leq 2
}B_{II_1I_2}\partial_{i_1}\left(Z^I_1(e^{-u}-1)\partial_{i_2}Z^{I_2}u\right)}+\vspace{3mm}\\
&  & {\displaystyle \sum_{|I_1|+|I_2|\leq |I|}\sum_{0\leq i_1,
i_2\leq 2
}\bar{B}_{II_1I_2}\partial_{i_1}Z^{I_1}u\partial_{i_2}Z^{I_2}u}+\vspace{3mm}\\
&  & {\displaystyle \sum_{|I_1|+|I_2|\leq |I|}\sum_{0\leq i_1,
i_2\leq 2
}\tilde{B}_{II_1I_2}\partial_{i_1}\left(Z^{I_1}(e^{-u}-1)\right)\partial_{i_2}Z^{I_2}u.}\end{array}$$
Let
\begin{equation}
Z^Iu=w_0+w_1+w_2+w_3,
\end{equation}
where $ w_0, w_1, w_2 $ and $ w_3 $ satisfy
\begin{equation}
\square w_0=0,\quad w_0|_{t=0}=Z^Iu(0,x),\quad \left.\frac{\partial
w_0}{\partial t}\right|_{t=0}=\frac{\partial (Z^Iu)}{\partial
t}(0,x),
\end{equation}
\begin{equation}
\square w_1=\sum_{|I_1|+|I_2|\leq |I|}\sum_{0\leq i_1, i_2\leq 2
}B_{II_1I_2}\partial_{i_1}\left(Z^I_1(e^{-u}-1)\partial_{i_2}Z^{I_2}u\right),\quad
w_1|_{t=0}=\left.\frac{\partial w_1}{\partial t}\right|_{t=0}=0,
\end{equation}
\begin{equation}
\square w_2=\sum_{|I_1|+|I_2|\leq |I|}\sum_{0\leq i_1, i_2\leq 2
}\bar{B}_{II_1I_2}\partial_{i_1}Z^{I_1}u\partial_{i_2}Z^{I_2}u,,\quad
w_2|_{t=0}=\left.\frac{\partial w_2}{\partial t}\right|_{t=0}=0,
\end{equation}
and
\begin{equation}
\square w_3=\sum_{|I_1|+|I_2|\leq |I|}\sum_{0\leq i_1, i_2\leq 2
}\tilde{B}_{II_1I_2}\partial_{i_1}\left(Z^{I_1}(e^{-u}-1)\right)\partial_{i_2}Z^{I_2}u,
\quad w_3|_{t=0}=\left.\frac{\partial w_3}{\partial
t}\right|_{t=0}=0,
\end{equation}
respectively. Thanks to Lemma 2.7 and (3.18), we have
\begin{equation} \|w_0(t, \cdot)\|_{L^2}\leq
\left\{\begin{array}{l} C\varepsilon\sqrt{\ln(2+t)} \quad (|x|\leq t),\\
C\varepsilon \quad (|x|\geq t).\end{array}\right.
\end{equation}
When $\varepsilon_0>0$ is suitably small, noting Lemmas 2.3, 2.5 and
using (3.19) and (3.2), we obtain
\begin{equation}\begin{array}{lll}
{\displaystyle
 \|w_1(t, \cdot)\|_{L^2}} & \leq & {\displaystyle C\left(\int^t_0\sum_{|I_1|+|I_2|\leq |I|}\sum_{0\leq i\leq 2
}\|Z^{I_1}(e^{-u}-1)\partial_iZ^{I_2}u(\tau,
\cdot)\|_{L^2}d\tau\right.+}\vspace{3mm}\\
&  &\qquad\quad {\displaystyle \left.\sum_{|I_1|+|I_2|\leq
|I|}\sum_{0\leq i \leq 2 }\|Z^{I_1}(e^{-u}-1)\partial_iZ^{I_2}u(0,
\cdot)\|_{L^2}\right)}\vspace{3mm}\\
& \leq & {\displaystyle
C\int^t_0\left[M_2(\tau)N_1(\tau)+M_1(\tau)N_2(\tau)\right]d\tau
+C\varepsilon.}\end{array}
\end{equation}
Noting Lemmas 2.9, 2.1 and using (3.20), (3.2) gives
\begin{equation}\begin{array}{lll}
{\displaystyle \|w_2(t, \cdot)\|_{L^2}} & \leq & {\displaystyle
C(1+t)^{\frac{1}{4}}\sum_{|I_1|+|I_2|\leq |I|,
|I_1|>|I_2|}\sum_{|J|\leq 1}\sum_{0\leq i_1, i_2\leq
2}\left\{\int^t_0(1+\tau)^{-\frac{1}{2}}\|Z^J\partial_{i_1}Z^{I_2}u(\tau,
\cdot)\|^2_{L^2}\right\}^{\frac{1}{2}}}\vspace{3mm}\\
& & \hskip 7cm {\displaystyle \times
\left\{\int^t_0\|\partial_{i_2}Z^{I_1}u(\tau,
\cdot)\|^2_{L^2}d\tau\right\}^{\frac{1}{2}}}+\vspace{3mm}\\
& & {\displaystyle C(1+t)^{\frac{1}{4}}\sum_{|I_1|+|I_2|\leq |I|,
|I_1|\leq|I_2|}\sum_{|J|\leq 1}\sum_{0\leq i_1, i_2\leq
2}\left\{\int^t_0(1+\tau)^{-\frac{1}{2}}\|Z^J\partial_{i_1}Z^{I_1}u(\tau,
\cdot)\|^2_{L^2}\right\}^{\frac{1}{2}}}\vspace{3mm}\\
& & \hskip 7cm {\displaystyle \times
\left\{\int^t_0\|\partial_{i_2}Z^{I_2}u(\tau,
\cdot)\|^2_{L^2}d\tau\right\}^{\frac{1}{2}}}\vspace{3mm}\\
& \leq &  {\displaystyle
C(1+t)^{\frac{1}{4}}\left\{\int^t_0(1+\tau)^{-\frac{1}{2}}M^2_1(\tau)d\tau\right\}^{\frac{1}{2}}
\left\{\int^t_0M^2_1(\tau)d\tau\right\}^{\frac{1}{2}}}.\end{array}\end{equation}
By Lemma 2.9, Minkowski inequality, Lemmas 2.1, 2.6 and (3.2), it
follows from (3.21) that
\begin{equation}\begin{array}{l}
{\displaystyle \|w_3(t, \cdot)\|_{L^2}\leq
C(1+t)^{\frac{1}{4}}\sum_{|I_1|+|I_2|\leq |I|,
|I_1|>|I_2|}\sum_{|J|\leq 1}\sum_{0\leq i_1, i_2\leq
2}\left\{\int^t_0(1+\tau)^{-\frac{1}{2}}\|Z^J\partial_{i_2}Z^{I_2}u(\tau,
\cdot)\|^2_{L^2}d\tau\right\}^{\frac{1}{2}}}\vspace{3mm}\\
\hskip 8cm {\displaystyle \times
\left\{\int^t_0\|\partial_{i_1}Z^{I_1}(e^{-u}-1)(\tau,
\cdot)\|_{L^2}d\tau\right\}^{\frac{1}{2}}}+\vspace{3mm}\\
\qquad\quad {\displaystyle C(1+t)^{\frac{1}{4}}\sum_{|I_1|+|I_2|\leq
|I|, |I_1|\leq|I_2|}\sum_{|J|\leq 1}\sum_{0\leq i_1, i_2\leq
2}\left\{\int^t_0(1+\tau)^{-\frac{1}{2}}\|Z^J\partial_{i_1}Z^{I_1}(e^{-u}-1)(\tau,
\cdot)\|^2_{L^2}d\tau\right\}^{\frac{1}{2}}}\vspace{3mm}\\
 \hskip 10cm {\displaystyle \times
\left\{\int^t_0\|\partial_{i_2}Z^{I_2}u(\tau,
\cdot)\|_{L^2}d\tau\right\}^{\frac{1}{2}}}\vspace{3mm}\\
\quad \leq  {\displaystyle C(1+t)^{\frac{1}{4}}\sum_{|I_1|+|I_2|\leq
|I|, |I_1|>|I_2|}\sum_{|J|\leq 1}\sum_{0\leq i_1, i_2\leq
2}\left\{\int^t_0(1+\tau)^{-\frac{1}{2}}\|Z^J\partial_{i_2}Z^{I_2}u(\tau,
\cdot)\|^2_{L^2}d\tau\right\}^{\frac{1}{2}}}\vspace{3mm}\\
 \hskip 5cm {\displaystyle \times
\left\{\int^t_0(\|Z^{I_1}((e^{-u}-1)\partial_{i_1}u)(\tau,
\cdot)\|^2_{L^2}+\|Z^{I_1}\partial_{i_1}u(\tau,
\cdot)\|^2_{L^2})d\tau\right\}^{\frac{1}{2}}}+\vspace{3mm}\\
\qquad\quad {\displaystyle C(1+t)^{\frac{1}{4}}\sum_{|I_1|+|I_2|\leq
|I|, |I_1|\leq|I_2|}\sum_{|J|\leq 1}\sum_{0\leq i_1, i_2\leq
2}\left\{\int^t_0(1+\tau)^{-\frac{1}{2}}\left[\|Z^JZ^{I_1}((e^{-u}-1)\partial_{i_1}u)(\tau,
\cdot)\|^2_{L^2}\right.\right.+}\vspace{3mm}\\
 \hskip 6cm {\displaystyle
\left.\left.\|Z^JZ^{I_1}\partial_{i_1}u(\tau,
\cdot)\|^2_{L^2}\right]d\tau\right\}^{\frac{1}{2}} \times
\left\{\int^t_0\|\partial_{i_2}Z^{I_2}u(\tau,
\cdot)\|_{L^2}d\tau\right\}^{\frac{1}{2}}}\vspace{3mm}\\
\quad \leq {\displaystyle
C(1+t)^{\frac{1}{4}}\left\{\int^t_0(1+\tau)^{-\frac{1}{2}}M_1^2(\tau)d\tau\right\}^{\frac{1}{2}}
\left\{\int^t_0\left[M_1^2(\tau)N^2_2(\tau)+M_2^2(\tau)N^2_1(\tau)+M_1^2(\tau)\right]d\tau\right\}^{\frac{1}{2}}+}\vspace{3mm}\\
\qquad\quad {\displaystyle
C(1+t)^{\frac{1}{4}}\left\{\int^t_0(1+\tau)^{-\frac{1}{2}}\left[M_1^2(\tau)N^2_2(\tau)+M_2^2(\tau)N^2_1(\tau)+M_1^2(\tau)
\right]d\tau\right\}^{\frac{1}{2}}\times
\left\{\int^t_0M_1^2(\tau)d\tau\right\}^{\frac{1}{2}}}\vspace{3mm}\\
\quad \leq  {\displaystyle
C(1+t)^{\frac{1}{4}}\left\{\int^t_0(1+\tau)^{-\frac{1}{2}}\left[M_1^2(\tau)N^2_2(\tau)+M_2^2(\tau)N^2_1(\tau)+
M_1^2(\tau)\right]d\tau\right\}^{\frac{1}{2}}}\vspace{3mm}\\
\hskip 4cm {\displaystyle \times
\left\{\int^t_0\left[M_1^2(\tau)N^2_2(\tau)+M_2^2(\tau)N^2_1(\tau)+M_1^2(\tau)\right]d\tau\right\}^{\frac{1}{2}}.}\end{array}\end{equation}
provided that $\varepsilon_0>0$ is suitably small. Thus, combining
(3.17) and (3.22)-(3.25) yields
\begin{equation}\begin{array}{l}
{\displaystyle M_2(t)\leq K_2\varepsilon \sqrt{\ln(2+t)}+
C\int^t_0(M_2(\tau)N_1(\tau)+M_1(\tau)N_2(\tau))d\tau + }\vspace{3mm}\\
\qquad\qquad {\displaystyle
C(1+t)^{\frac{1}{4}}\left\{\int^t_0(1+\tau)^{-\frac{1}{2}}\left[M_1^2(\tau)N^2_2(\tau)+M_2^2(\tau)
N^2_1(\tau)+M_1^2(\tau)\right]d\tau\right\}^{\frac{1}{2}}}\vspace{3mm}\\
\hskip 4cm {\displaystyle
\times\left\{\int^t_0\left[M_1^2(\tau)N^2_2(\tau)+M_2^2(\tau)N^2_1(\tau)+M_1^2(\tau)\right]
d\tau\right\}^{\frac{1}{2}}.}\end{array}\end{equation}

We now estimate $N_2(t)$.

Using Lemma 2.7, we obtain from (3.18) that
\begin{equation}
(1+t)^{\frac{1}{2}}\|w_0(t, \cdot)\|_{L^\infty}\leq C\varepsilon.
\end{equation}
Noting Lemmas 2.3, 2.5 and using (3.19), (3.2), when
$\varepsilon_0>0$ is suitably small,  we have
\begin{equation}\begin{array}{l}
{\displaystyle (1+t)^{\frac{1}{2}}\|w_1(t, \cdot)\|_{L^\infty}\leq
C\int^t_0(1+\tau)^{\frac{1}{2}}\sum_{|I_1|+|I_2|\leq|I|}\sum_{0\leq
i\leq 2}\|Z^{I_1}(e^{-u}-1)\partial_iZ^{I_2}u(\tau,
\cdot)\|_{L^\infty}d\tau}+\vspace{3mm}\\
\qquad\qquad {\displaystyle
C\int^t_0(1+\tau)^{-\frac{3}{2}}\sum_{|I_1|+|I_2|\leq|I|}\sum_{|J|\leq
3}\sum_{0\leq i\leq
2}\|Z^J\left(Z^{I_1}(e^{-u}-1)\partial_iZ^{I_2}u\right)(\tau,
\cdot)\|_{L^1}d\tau}\vspace{3mm}\\
\qquad\quad {\displaystyle\leq
C\left\{\int^t_0(1+\tau)^{\frac{1}{2}}N_1(\tau)N_2(\tau)d\tau+\int^t_0(1+\tau)^{-\frac{3}{2}}M_1(\tau)M_2(t)d\tau\right\}.
}\end{array}\end{equation} Noting Lemma 2.9 and using (3.20) and
(3.2), we obtain
\begin{equation}\begin{array}{lll}
(1+t)^{\frac{1}{2}}\|w_2(t, \cdot)\|_{L^\infty} & \leq &
{\displaystyle C \sum_{|I_1|+|I_2|\leq|I|}\sum_{|J|=1}\sum_{0\leq
i_1, i_2\leq 2}\left\{\int^t_0\|Z^J(\partial_{i_1}Z^{I_1}u)(\tau,
\cdot)\|^2_{L^2}\frac{d\tau}{\sqrt{1+\tau}}\right\}^{\frac{1}{2}}}\vspace{3mm}\\
& & {\displaystyle \times
\left\{\int^t_0\|Z^J(\partial_{i_2}Z^{I_2}u)(\tau,
\cdot)\|^2_{L^2}\frac{d\tau}{\sqrt{1+\tau}}\right\}^{\frac{1}{2}}}\vspace{3mm}\\
& \leq & {\displaystyle
C\int^t_0M^2_1(\tau)\frac{d\tau}{\sqrt{1+\tau}}.}\end{array}\end{equation}
By Lemmas 2.9, 2.1, 2.6, Minkowski inequality and (3.2), when
$\varepsilon_0>0$ is suitably small, we obtain from (3.21) that
\begin{equation}\begin{array}{l}
{\displaystyle (1+t)^{\frac{1}{2}}\|w_3(t, \cdot)\|_{L^\infty}\leq
C\sum_{|I_1|+|I_2|\leq|I|}\sum_{|J|\leq 1}\sum_{0\leq i_1, i_2\leq
2}\left\{\int^t_0\|Z^J\left(\partial_{i_1}Z^{I_1}(e^{-u}-1)\right)(\tau,
\cdot)\|^2_{L^2}\frac{d\tau}{\sqrt{1+\tau}}\right\}^{\frac{1}{2}}}\vspace{3mm}\\
\hskip 7cm {\displaystyle \times
\left\{\int^t_0\|Z^J(\partial_{i_2}Z^{I_2}u)(\tau,
\cdot)\|^2_{L^2}\frac{d\tau}{\sqrt{1+\tau}}\right\}^{\frac{1}{2}}}\vspace{3mm}\\
\quad \leq {\displaystyle
C\left\{\int^t_0(1+\tau)^{-\frac{1}{2}}\left[M_1^2(\tau)N^2_2(\tau)+M_2^2(\tau)N^2_1(\tau)+
M_1^2(\tau)\right]d\tau\right\}^{\frac{1}{2}}
\left\{\int^t_0(1+\tau)^{-\frac{1}{2}}M_1^2(\tau)d\tau\right\}^{\frac{1}{2}}.}\end{array}\end{equation}
Collecting (3.17) and (3.27)-(3.30) gives
\begin{equation}\begin{array}{l}
{\displaystyle (1+t)^{\frac{1}{2}}N_2(t)\leq K_3\varepsilon+
C\int^t_0(1+\tau)^{\frac{1}{2}}N_1(\tau)N_2(\tau)d\tau+C\int^t_0(1+\tau)^{-\frac{3}{2}}M_1(\tau)M_2(t)d\tau
}\vspace{3mm}\\
 {\displaystyle \qquad\quad
+C\left\{\int^t_0(1+\tau)^{-\frac{1}{2}}(M_1^2(\tau)N^2_2(\tau)+M_2^2(\tau)N^2_1(\tau)+M_1^2(\tau))d\tau\right\}^{\frac{1}{2}}
\left\{\int^t_0(1+\tau)^{-\frac{1}{2}}M_1^2(\tau)d\tau\right\}^{\frac{1}{2}}.}\end{array}\end{equation}

Since, for the time being it supposed that (3.3) holds, (3.12) and
(3.13) are true, provided that $\varepsilon_0$ is suitably small,
and then, by (3.15), (3.26) and (3.31) it holds that
\begin{equation}
M_1(t)\leq
K_1\varepsilon+C(M_2N_1+M_1N_2+M_1N_1)\varepsilon^2(1+t)^{\frac{3}{4}},
\end{equation}
\begin{equation}
M_2(t)\leq
K_2\varepsilon\sqrt{\ln(2+t)}+C(M_2N_1+M_1N_2+M_1N_1)\varepsilon^2(1+t)
\end{equation}
and
\begin{equation}
(1+t)^{\frac{1}{2}}N_2(t)\leq
K_3\varepsilon+C(N_1N_2+M_1M_2+M_2N_1+M_1N_2+M^2_1)\varepsilon^2\sqrt{1+t}.
\end{equation}
On the other hand, by Lemma 2.10, we have
\begin{equation}
(1+t)^{\frac{1}{2}}N_{1}(t)\leq CM_1(t).
\end{equation}
Thus, choosing
$$M_1\geq 4K_1, \quad M_2\geq 4K_2, \quad N_2\geq 4K_3, N_1\geq 2CM_1,$$
we obtain from (3.32)-(3.35)
that
\begin{equation}\left\{\begin{array}{l}
M_1(t)\leq \frac{1}{2}M_1\varepsilon,\vspace{2mm}\\
M_2(t)\leq\frac{1}{2}M_2\varepsilon (1+t)^{\frac{1}{4}},\vspace{2mm}\\
(1+t)^{\frac{1}{2}}N_1(t)\leq
\frac{1}{2}N_1\varepsilon,\vspace{2mm}\\
(1+t)^{\frac{1}{2}}N_2(t)\leq
\frac{1}{2}N_2\varepsilon,\end{array}\right.\end{equation} provided
that $\varepsilon_0$ and $\mu $ are suitably small and satisfy
$$\varepsilon_0(1+t)^{\frac{3}{4}}\leq \mu.$$
 Take $\delta=\mu^{\frac{4}{3}}$. Thus, the proof of Theorem 3.1
is completed. $\quad\quad\quad\blacksquare$ \vskip 4mm

\noindent{\Large \textbf{Acknowledgements.}} The authors thank Dr.
Zhen Lei and Prof. Yi Zhou for their valuable suggestions. This work
was supported in part by the NNSF of China (Grant No. 10971190) and
the Qiu-Shi Professor Fellowship from Zhejiang University,
China.%; the work of Liu was supported in part by the NSF and NSF of China.

\end{document}